\documentclass[11pt,draft]{article}
\usepackage[a4paper,centering,vscale=0.75,hscale=0.7]{geometry}
\usepackage{mathtools}
\usepackage{amsthm,amssymb}
\usepackage{mathrsfs}

\newtheorem{thm}{Theorem}[section]
\newtheorem{defi}[thm]{Definition}
\theoremstyle{remark}
\newtheorem{rmk}[thm]{Remark}

\renewcommand{\leq}{\leqslant}
\renewcommand{\geq}{\geqslant}
\newcommand{\dom}{\mathsf{D}}
\newcommand{\E}{\mathop{{}\mathbb{E}}}
\newcommand{\cF}{\mathscr{F}}
\newcommand{\cL}{\mathscr{L}}
\renewcommand{\P}{\mathbb{P}}
\newcommand{\erre}{\mathbb{R}}

\newcommand{\embed}{\hookrightarrow}

\DeclarePairedDelimiter\norm{\lVert}{\rVert}
\DeclarePairedDelimiterX\ip[2]{\langle}{\rangle}{#1,#2}
\DeclarePairedDelimiterX\ipp[2]{\langle\!\langle}{\rangle\!\rangle}{#1,#2}

\numberwithin{equation}{section}

\title{An alternative proof of well-posedness of\\
  stochastic evolution equations in the variational setting}
\author{Carlo Marinelli\thanks{Department of Mathematics, University
    College London, Gower Street, London WC1E 6BT, United
    Kingdom. URL: \texttt{http://goo.gl/4GKJP}} \and Luca
  Scarpa\thanks{Faculty of Mathematics, University of Vienna,
    Oskar-Morgenstern-Platz 1, 1090 Vienna, Austria. E-mail:
    \texttt{luca.scarpa@univie.ac.at}} \and Ulisse
  Stefanelli\thanks{Faculty of Mathematics, University of Vienna,
    Oskar-Morgenstern-Platz 1, A-1090 Vienna, Austria; Vienna Research
    Platform on Accelerating Photoreaction Discovery, University of
    Vienna, W\"ahringer Str. 17, 1090 Vienna, Austria; Istituto di
    Matematica Applicata e Tecnologie Informatiche ``E. Magenes'' --
    CNR, via Ferrata 1, I-27100 Pavia, Italy.  E-mail:
    \texttt{ulisse.stefanelli@univie.ac.at}}}
\date{September 21, 2020}

\begin{document}
\maketitle

\begin{abstract}
  We present a new proof of well-posedness of stochastic evolution
  equations in variational form, relying solely on a (nonlinear)
  infinite-dimensional approximation procedure rather than on
  classical finite-dimensional projection arguments of Galerkin type.
\end{abstract}


\section{Introduction}
\label{sec:intro}
Let us consider a stochastic evolution equation of the type
\begin{equation}
  \label{eq:1}
  dX(t) + A(t,X(t))\,dt = B(t,X(t))\,dW(t), \qquad X(0)=X_0,
\end{equation}
in the so-called variational setting, i.e. where $A$ is a random
time-dependent nonlinear maximal monotone operator from a reflexive
Banach space $V$ to its dual $V'$, with $V$ densely and continuously
embedded in a Hilbert space $H$. Moreover, $W$ is a cylindrical Wiener
process (possibly defined on a further separable Hilbert space), and
$B$ is a random time-dependent map with values in a suitable space of
Hilbert-Schmidt operators. Precise assumptions on the data of the
problem are given in \S\ref{sec:main} below.

This class of equations was introduced and studied by Pardoux in
\cite{Pard}, extending to the stochastic setting the classical
well-posedness results by Lions (see, e.g., \cite{Lions:q}) for
equations without noise. More precisely, in \cite{Pard} the operator
$A$ is time-dependent but non-random and the pair $(A,B)$ needs to
satisfy coercivity and boundedness assumptions, complemented by a
local Lipschitz continuity condition on $B$. The general case where
$A$ can be random was considered by Krylov and Rozovski\u{\i} in
\cite{KR-spde}, who also showed that the local Lipschitz continuity on
$B$ is not needed (see also \cite{LiuRo,RozLot} for comprehensive
treatments and recent developments). Under coercivity and boundedness
assumptions on $A$ and a Lipschitz continuity assumption on $B$ (a set
of hypotheses to which we shall refer to as \emph{disjoint}
assumptions), Pardoux \cite[Chapter~3, {\S}1]{Pard} proved
well-posedness of \eqref{eq:1} by a clever, but at the same time
natural extension of the deterministic theory, employing an
infinite-dimensional argument based on Picard iterations in suitable
spaces of processes. As a second step, assuming that the pair $(A,B)$
satisfies a \emph{joint} coercivity and boundedness assumption and
that $B$ is locally Lipschitz continuous, well-posedness for
\eqref{eq:1} is proved by a different method, i.e. by
finite-dimensional approximations of Galerkin type (see
\cite[Chapter~3, {\S}3]{Pard}). As mentioned above, the joint
assumption is shown to imply well-posedness without any local
Lipschitz continuity condition in $B$ in \cite{KR-spde}, again using
finite-dimensional approximations.

Our goal is to show that existence of solutions under the general
joint assumptions on $(A,B)$, as in \cite{KR-spde}, can be obtained
relying only on infinite-dimensional arguments, i.e. in the same
spirit of the approach adopted in the first part of
\cite{Pard}. Moreover, since the existence proof in \cite{KR-spde}
relies on quite advanced results for finite-dimensional stochastic
differential equations, the alternative proof provided here could also
be seen as a {\it simpler} proof.
The main idea is, roughly speaking, to regularize the operator $B$
through the resolvent of the operator $A$. The corresponding
regularized problem is then shown to satisfy the stronger
\emph{disjoint} hypotheses, so that it can be solved, as in
\cite[Chapter~3, \S~1]{Pard}, using infinite-dimensional techniques
only. Uniform estimates on the solutions to the regularized equations
are then established, which allow to pass to the limit obtaining a
solution to the original problem.

Even though the use of Yosida approximations is a standard tool in the
field of nonlinear deterministic and stochastic equations (see, e.g.,
\cite{Barbu:type,Bmax,Govi,LiuSte} for just a few examples among an
enormous literature), it seems that the type of approximation
introduced here is not found elsewhere. Let us also mention that
another approach to stochastic equations in variational form, namely
by reduction to the deterministic case, is developed in
\cite[{\S}4.4]{Barbu:type}, where, however, only the case of additive
noise is considered. Moreover, well-posedness for certain classes of
stochastic equations with multi-valued nonlinear drift term is
obtained in \cite{LiuSte} using the results in \cite{KR-spde,Pard} as
starting point. The main point in \cite{LiuSte} is to regularize the
multi-valued drift coefficient by its Yosida approximation, thus
obtaining a family of well-posed equations with single-valued drift
satisfying the assumptions of the classical variational framework, and
to show that the corresponding approximate solutions converge, under
appropriate assumptions, among which the Lipschitz continuity of the
diffusion coefficient, to a process solving the original equation.


\section{Setting and main results}
\label{sec:main}
Throughout the paper, $(\Omega,\cF,(\cF_t)_{t\in[0,T]}, \P)$ stands
for a filtered probability space satisfying the so-called {\it usual}
conditions, where $T>0$ is a fixed final time, on which all random
elements will be defined.  Equality of processes is always meant in
the sense of indistinguishability, unless otherwise stated. Moreover,
$U$ is a separable Hilbert space and $W$ is a cylindrical Wiener
process on it; $H$ is a separable Hilbert space identified with its dual,
and $V$ is a separable reflexive Banach space continuously and densely
embedded in $H$, so that, denoting the (topological) dual of $V$ by
$V'$, $V \embed H \embed V'$ is a Gelfand triple. The scalar product
and norm of $H$ will be denoted by $\ip{\cdot}{\cdot}$ and
$\norm{\cdot}$, respectively, while the norms of all other Banach
spaces will be indicated by subscripts. Since the duality form between
$V$ and $V'$ agrees with the scalar product of $H$ in the usual sense,
we shall denote the former by $\ip{\cdot}{\cdot}$ as well.  If $E_1$
and $E_2$ are Hilbert spaces, the space of Hilbert-Schmidt operators
from $E_1$ to $E_2$ will be denoted by $\cL^2(E_1;E_2)$.

\medskip

The following assumptions will be in force throughout the paper.
\begin{description}
\item[(I)] The operator $A:\Omega\times[0,T]\times V\to V'$ is
  progressively measurable and hemicontinuous, i.e., for every
  $x\in V$ the $V'$-valued process $A(\cdot,\cdot,x)$ is progressively
  measurable and the map
  \[
  \erre \ni r \longmapsto \ip{A(\omega,t,x+ry)}{z}
  \]
  is continuous for every $\omega \in \Omega$, $t \in [0,T]$, and
  $x,y,z \in V$.
\item[(II)] The operator $B\colon \Omega \times [0,T] \times V \to \cL^2(U;H)$
  is progressively measurable, i.e., for every $x \in V$ the
  $\cL^2(U;H)$-valued process $B(\cdot,\cdot,x)$ is progressively
  measurable.
\item[(III)] There exist constants $c_1>0$, $c_2\geq0$,
  $p \in \mathopen]1,+\infty\mathclose[$ and an adapted process
  $f \in L^1(\Omega\times(0,T))$ such that
  \begin{gather*}
    \ip[\big]{A(\omega,t,x)-A(\omega,t,y)}{x-y}
    - \frac12 \norm[\big]{B(\omega,t,x)-B(\omega,t,y)}_{\cL^2(U;H)}^2
    \geq - c_2\norm{x-y}^2,\\
    \ip[\big]{A(\omega,t,x)}{x}
    - \frac12 \norm[\big]{B(\omega,t,x)}_{\cL^2(U;H)}^2
    \geq c_1\norm{x}_V^p - c_2\norm{x}^2 - f(\omega,t)
  \end{gather*}
  for every $x,y\in V$, $t \in [0,T]$ and $\omega \in \Omega$.
\item[(IV)] There exist a constant $C>0$ and an adapted process
  $g \in L^1(\Omega \times (0,T))$ such that, setting
  $q:=p/(p-1)$,
  \[
  \norm[\big]{A(\omega,t,x)}_{V'}^q \leq C\norm{x}_V^p + g(\omega,t)
  \]
  for every $x \in V$, $t \in [0,T]$ and $\omega \in \Omega$.
  \item[(V)] $X_0 \in L^2(\Omega,\cF_0;H)$.
\end{description}

We can give now the definition of strong solution for the equation.
\begin{defi}
  A strong solution to \eqref{eq:1} is a $V$-valued progressively
  measurable process $X$ such that
  \begin{gather*}
  X \in L^0(\Omega;C([0,T];H)) \cap L^0(\Omega;L^p(0,T; V)),\\
  B(\cdot,\cdot,X)\in L^0(\Omega;L^2(0,T;\cL^2(U;H))),
  \end{gather*}
  and
  \[
    X + \int_0^\cdot A(s,X(s))\,ds = X_0 + \int_0^\cdot
    B(s,X(s))\,dW(s),
  \]
  as an identity in the sense of indistinguishable $V'$-valued
  processes.
\end{defi}

The classical well-posedness result \cite{KR-spde, Pard}
for equation \eqref{eq:1} is as follows.
\begin{thm}
  \label{th:wp}
  There exists a unique strong solution $X$ to \eqref{eq:1}. Moreover,
  \[
  X \in L^2(\Omega;C([0,T];H)) \cap L^p(\Omega;L^p(0,T; V))
  \]
  and the solution map 
  \begin{align*}
    L^2(\Omega,\cF_0; H)
    &\longrightarrow C([0,T];L^2(\Omega;H))\\
    X_0 &\longmapsto X
  \end{align*}
  is Lipschitz continuous.
\end{thm}

As discussed above, in the next section we show that
Theorem~\ref{th:wp} can be proved relying only on infinite-dimensional
arguments.
\begin{rmk}
  The proof of uniqueness of strong solution crucially relies on an
  It\^o formula for the square of the $H$-norm. In \cite{Pard} such
  formula is obtained under the assumption that an operator
  $C\colon V \to V'$ exists satisfying monotonicity, coercivity and
  boundedness conditions (see~\cite[p.~57]{Pard}).  These conditions
  coincide with those assumed here on $A(\omega,t)$ for all
  $(\omega,t) \in \Omega \times [0,T]$, hence are automatically
  verified. In general, the formula remains valid, without any
  connection to a specific stochastic equation, if the duality map
  $J \colon V \to V'$ is single-valued, which is the case if, e.g.,
  $V'$ is a strictly convex Banach space. Such an assumption is always
  satisfied in all applications to SPDEs we know of. On the other
  hand, the It\^o formula in \cite{KR-spde} does not require any
  ``geometric'' assumption on the Banach space $V$, but its proof is
  rather involved and relies on finite-dimensional projections. A
  simple proof of It\^o's formula for the square of the $H$-norm in
  the variational setting, which relies just on infinite-dimensional
  arguments, is available, to the best of our knowledge, only in the
  case where $V$ is a Hilbert space (see~\cite{Kry:shortIto}, as well
  as \cite{cm:semimg}).
\end{rmk}


\section{Proof of Theorem~\ref{th:wp}}
\label{sec:proof1}
With $\omega \in \Omega$ and $t \in [0,T]$ arbitrary but fixed,
assumptions {\bf (I)} and {\bf (III)} imply that
$\tilde{A} := A(\omega,t,\cdot) + c_2 I\colon V \to V'$ is maximal
monotone.  Let us show that the part of $\tilde{A}$ in $H$, denoted by
$\tilde{A}_H$, is maximal monotone (as an operator in $H$) with domain
$\dom(\tilde{A}) := \{x \in V: A(\omega,t,x) \in H\}$. It suffices to
show that, for any $y \in H$, the equation
\[
  x + \tilde{A}x = y
\]
admits a solution $x \in \dom(\tilde{A})$. Since $\tilde{A}$ is
coercive on $V$ by assumption {\bf (III)}, it follows by maximal
monotonicity that the equation admits a (unique) solution $x \in
V$. This obviously implies $\tilde{A}x \in H$, i.e.
$x \in \dom(\tilde{A})$. We have hence shown that $\tilde{A}_H$ is a
maximal monotone operator on $H$.

For every $\lambda>0$ we define the resolvent operator
\[
  J_\lambda \colon \Omega \times [0,T] \times H \longrightarrow V
\]
of $\tilde{A}$ setting, for every $x \in H$, $t \in [0,T]$ and
$\omega \in \Omega$,
\[
  J_\lambda(\omega,t,x) + \lambda\tilde
  A(\omega,t,J_\lambda(\omega,t,x)) = x.
\]
The maximal monotonicity of $\tilde{A}_H$ implies that, for every
$(\omega,t) \in \Omega \times [0,T]$, $J_\lambda(\omega,t,\cdot)$ is a
contraction and converges pointwise to the identity map of $H$ as
$\lambda \to 0$.

Moreover, we define the Yosida approximation of $\tilde{A}$ as the map
\begin{align*}
  \tilde{A}_\lambda \colon \Omega \times[0,T] \times H
  &\longrightarrow H,\\
  (\omega,t,x)
  &\longmapsto \tilde{A}(\omega,t,J_\lambda(\omega,t,x)).
\end{align*}
It follows by the contraction property of $J_\lambda$ that, for every
$(\omega,t) \in \Omega \times [0,T]$,
$\tilde{A}_\lambda(\omega,t,\cdot)$ is Lipschitz continuous with
Lipschitz constant bounded by $1/\lambda$.

\subsection{Regularized equation}
Let us introduce the family of operators indexed by $\lambda>0$
\begin{align*}
  B_\lambda \colon \Omega \times[0,T] \times H
  &\longrightarrow \cL^2(U;H)\\
  (\omega,t,x) &\longmapsto B(\omega,t,J_\lambda(\omega,t,x)).
\end{align*}
Since $A$, hence also $J_\lambda$, and $B$ are progressively
measurable, $B_\lambda$ is progressively measurable as well for every
$\lambda>0$. Moreover, $B_\lambda$ is Lipschitz continuous in its
third argument, uniformly with respect to the other ones: in fact,
thanks to assumption {\bf (III)}, one has
\begin{align*}
  &\frac12 \norm[\big]{B_\lambda(\omega,t,x)%
    - B_\lambda(\omega,t,y)}_{\cL^2(U;H)}^2\\
  &\hspace{3em} =
    \frac12 \norm[\big]{B(\omega,t,J_\lambda(\omega,t,x))%
    - B(\omega,t,J_\lambda(\omega,t,y))}_{\cL^2(U;H)}^2\\
  &\hspace{3em} \leq \ip[\big]{A(\omega,t,J_\lambda(\omega,t,x))%
    - A(\omega,t,J_\lambda(\omega,t,y))}%
    {J_\lambda(\omega,t,x)-J_\lambda(\omega,t,y)}\\
  &\hspace{5em} + c_2\norm{J_\lambda(\omega,t,x) - J_\lambda(\omega,t,y)}^2\\
  &\hspace{3em} = \ip[\big]{\tilde{A}_\lambda(\omega,t,x)%
    - \tilde{A}_\lambda(\omega,t,y)}{x-y}%
    \leq \frac{1}{\lambda} \norm{x-y}^2
\end{align*}
for every $\omega \in \Omega$, $t \in [0,T]$, and $x,y \in H$.

Therefore, since $\tilde{A}_\lambda$ is also a Lipschitz continuous
operator on $H$, well-posedness results for stochastic differential
equations on Hilbert spaces (see, e.g., \cite[{\S}34]{Met}) yields the
existence (and uniqueness) of a predictable process
$X_\lambda \in L^2(\Omega;C([0,T];H))$ such that
\begin{equation}
  \label{eq:1app}
  X_\lambda + \int_0^\cdot \tilde{A}_\lambda(s,X_\lambda(s))\,ds 
  = X_0 + c_2 \int_0^\cdot X_\lambda(s)\,ds
  + \int_0^\cdot B_\lambda(s,X_\lambda(s))\,dW(s).
\end{equation}

\subsection{A priori estimates}
We are going to prove estimates on the family of solutions
$(X_\lambda)$ to the regularized equations obtained above that are
uniform with respect to $\lambda$.

The integration by parts formula for Hilbert space valued
semimartingales yields
\begin{align*}
  &\frac12 \norm{X_\lambda}^2
    + \int_0^\cdot \ip[\big]{\tilde{A}_\lambda(s,X_\lambda(s))}{X_\lambda(s)}\,ds
    - \frac12 \int_0^\cdot \norm[\big]{B_\lambda(s,X_\lambda(s))}^2_{\cL^2(U;H)}\,ds\\
  &\hspace{3em}= \frac12 \norm{X_0}^2
    + c_2 \int_0^\cdot \norm{X_\lambda(s)}^2\,ds
    + \int_0^\cdot X_\lambda(s) B_\lambda(s,X_\lambda(s))\,dW(s),
\end{align*}
where, in the stochastic integral, $X_\lambda$ is treated as a process
taking values in the dual of $H$.

From now on we shall occasionally suppress the explicit indication of
the dependence on $\omega$ and $t$ for processes and operators for
notational simplicity.  Note that, by definition of Yosida
approximation,
\begin{equation}
  \label{eq:yo}
\begin{split}
  \ip[\big]{\tilde{A}_\lambda x}{x} =
  \ip[\big]{\tilde{A}J_\lambda x}{x}
  &= \ip[\big]{\tilde{A}J_\lambda x}{J_\lambda x}
  + \ip[\big]{\tilde{A}J_\lambda x}{x - J_\lambda x}\\
  &= \ip[\big]{\tilde{A}J_\lambda x}{J_\lambda x}
  + \lambda \norm[\big]{\tilde{A}_\lambda x}^2
\end{split}
\end{equation}
for every $x \in H$. Therefore, denoting the norm of $\cL^2(U;H)$ by
$\norm{\cdot}_2$ for brevity, one has, thanks to assumption {\bf (III)},
\begin{align*}
  \ip[\big]{\tilde{A}_\lambda x}{x} - \frac12 \norm[\big]{B_\lambda(x)}_2^2
  &= \ip[\big]{\tilde{A}J_\lambda x}{J_\lambda x}
  - \frac12 \norm[\big]{B(J_\lambda x)}_2^2
    + \lambda \norm[\big]{\tilde{A}_\lambda x}^2\\
  &\geq c_1 \norm[\big]{J_\lambda x}_V^p
    + \lambda \norm[\big]{\tilde{A}_\lambda x}^2 - f,
\end{align*}
which in turn yields
\begin{align*}
  &\frac12 \norm{X_\lambda}^2
    + c_1 \int_0^\cdot \norm[\big]{J_\lambda X_\lambda(s)}_V^p\,ds
  + \lambda \int_0^\cdot \norm[\big]{\tilde{A}_\lambda X_\lambda(s)}^2\,ds\\
  &\hspace{3em} \leq \frac12 \norm{X_0}^2
    + \int_0^\cdot f(s)\,ds + c_2 \int_0^\cdot \norm{X_\lambda(s)}^2\,ds\\
  &\hspace{5em} + \int_0^\cdot X_\lambda(s) B_\lambda(X_\lambda(s))\,dW(s),
\end{align*}
where the stochastic integral on the right-hand side is a martingale
because $B_\lambda$ is Lipschitz continuous.  In particular, taking
expectation on both sides,
\begin{equation}
  \label{eq:ittica}
\begin{split}
  &\frac12 \E\norm{X_\lambda(t)}^2
    + c_1 \E\int_0^t \norm[\big]{J_\lambda X_\lambda(s)}_V^p\,ds
  + \lambda \E\int_0^t \norm[\big]{\tilde{A}_\lambda X_\lambda(s)}^2\,ds\\
  &\hspace{3em} \leq \frac12 \E\norm{X_0}^2
    + \E\int_0^t f(s)\,ds + c_2 \E\int_0^t \norm{X_\lambda(s)}^2\,ds
\end{split}
\end{equation}
for all $t \in [0,T]$. This implies that, for any interval
$[0,T_0] \subseteq [0,T]$,
\[
  \sup_{t \leq T_0} \E\norm{X_\lambda(t)}^2
  \leq \E\norm{X_0}^2 + 2\E\int_0^T f(s)\,ds
  + 2c_2 \E\int_0^{T_0} \Bigl(\sup_{r \leq s} \E\norm{X_\lambda(r)}^2\Bigr)\,ds,
\]
so that, by Gronwall's inequality, $(X_\lambda)$ is bounded in
$C([0,T];L^2(\Omega;H))$. From this and \eqref{eq:ittica} it
immediately follows that $(\lambda^{1/2} \tilde{A}_\lambda X_\lambda)$
is bounded in $L^2(\Omega;L^2(0,T;H))$ and that
$(J_\lambda X_\lambda)$ is bounded in $L^p(\Omega;L^p(0,T;V))$. The
latter implies, thanks to assumption {\bf (IV)}, that
$(\tilde{A}_\lambda X_\lambda)$ is bounded in
$L^q(\Omega;L^q(0,T;V'))$. Moreover, since, by assumption
{\bf (III)},
\[
  \frac12 \norm[\big]{B_\lambda(x)}_2^2
  = \frac12 \norm[\big]{B(J_\lambda x)}^2_2
  \leq \ip[\big]{\tilde{A}J_\lambda x}{J_\lambda x} + f
  = \ip[\big]{\tilde{A}_\lambda x}{J_\lambda x} + f,
\]
H\"older's inequality yields
\begin{align*}
  \E\int_0^T \norm[\big]{B_\lambda(X_\lambda(s))}_2^2\,ds
  &\leq \E\int_0^T \norm[\big]{\tilde{A}_\lambda X_\lambda}_{V'}%
    \norm[\big]{J_\lambda X_\lambda}_V
    + \norm[\big]{f}_{L^1(\Omega \times [0,T])}\\
  &\leq \norm[\big]{\tilde{A}_\lambda X_\lambda}_{L^q(\Omega;L^q(0,T;V'))}
    \norm[\big]{J_\lambda X_\lambda}_{L^p(\Omega;L^p(0,T;V))}
    + \norm[\big]{f}_{L^1(\Omega \times [0,T])},
\end{align*}
thus proving that $(B_\lambda(X_\lambda))$ is bounded in
$L^2(\Omega;L^2(0,T;\cL^2(U;H)))$.

\subsection{Construction of a solution and its uniqueness}
The boundedness of various families of processes indexed by $\lambda$
obtained above imply, by well-known compactness properties in weak and
weak* topologies, that there exist measurable and adapted processes
\begin{align*}
  X &\in L^\infty(0,T;L^2(\Omega;H)),\\
  \bar{X} &\in L^p(\Omega;L^p(0,T;V)),\\
  Y &\in L^q(\Omega;L^q(0,T;V')),\\
  G &\in L^2(\Omega;L^2(0,T;\cL^2(U;H)))
\end{align*}
such that
\begin{alignat*}{2}
  X_\lambda &\longrightarrow X &\qquad
  &\text{weakly* in } L^\infty(0,T;L^2(\Omega;H)),\\
  J_\lambda X_\lambda &\longrightarrow \bar{X}
  &&\text{weakly in } L^p(\Omega; L^p(0,T; V)),\\
  \lambda \tilde{A}_\lambda X_\lambda &\longrightarrow 0
  &&\text{in } L^2(\Omega;L^2(0,T;H)),\\
  \tilde{A}_\lambda X_\lambda &\longrightarrow Y
  &&\text{weakly in } L^q(\Omega; L^q(0,T; V')),\\
  B_\lambda(\cdot,X_\lambda) &\longrightarrow G
  &&\text{weakly in } L^2(\Omega; L^2(0,T; \cL^2(U;H))).
\end{alignat*}
Since, by definition of Yosida approximation,
$X_\lambda - J_\lambda X_\lambda = \lambda \tilde{A}_\lambda
X_\lambda$, the boundedness of
$(\lambda^{1/2}\tilde{A}_\lambda X_\lambda)$ in
$L^2(\Omega;L^2(0,T;H))$ implies that
\[
  X_\lambda - J_\lambda X_\lambda \longrightarrow 0 \qquad
  \text{in } L^2(\Omega;L^2(0,T;H)),
\]
hence $X$ and $\bar{X}$ are equal
$\mathbb{P} \otimes \operatorname{Leb}$-a.e. in $\Omega \times [0,T]$
and belong to
$L^\infty(0,T;L^2(\Omega;H)) \cap L^p(\Omega;L^p(0,T;V))$.
Recalling that the linear operator
\begin{align*}
  L^r(\Omega;L^1(0,T;E)) &\longrightarrow L^r(\Omega;C([0,T];E))\\
  w &\longmapsto \int_0^\cdot w(s)\,ds
\end{align*}
is continuous for any $r \in \mathopen[1,\infty\mathclose[$ and any
Banach space $E$, hence also weakly continuous, it follows from the
above convergence results that
\begin{alignat*}{2}
  \int_0^\cdot X_\lambda(s)\,ds &\longrightarrow 
  \int_0^\cdot X(s)\,ds
  &\qquad &\text{weakly in } L^p(\Omega; C([0,T]; V)),\\
  \int_0^\cdot\tilde{A}_\lambda X_\lambda(s)\,ds &\longrightarrow 
  \int_0^\cdot Y(s)\,ds
  &&\text{weakly in } L^q(\Omega; C([0,T]; V')).
\end{alignat*}
Similarly, since the stochastic integral operator
\begin{align*}
  L^2(\Omega;L^2(0,T;\cL^2(U;H)) &\longrightarrow L^2(\Omega; C([0,T]; H))\\
  C &\longmapsto \int_0^\cdot C(s)\,dW(s)
\end{align*}
is linear and continuous, hence weakly continuous, one gets that
\[
  \int_0^\cdot B_\lambda(s,X_\lambda(s))\,dW(s) \longrightarrow 
  \int_0^\cdot G(s)\,dW(s)
  \qquad \text{weakly in } L^2(\Omega; C([0,T]; H)).
\]
This implies, taking the weak limit in the space
$L^q(\Omega; C([0,T]; V'))$ as $\lambda \to 0$ in equation
\eqref{eq:1app},
\begin{equation}
  \label{eq:lim}
  X + \int_0^\cdot Y(s)\,ds = X_0 + c_2\int_0^\cdot X(s)\,ds
  + \int_0^\cdot G(s)\,dW(s) \qquad\text{in } V'.
\end{equation}
We are going to show, using an adaptation of a classical argument from
the theory of maximal monotone operators (cf.~\cite{KR-spde} as well
as \cite[\S~4.2]{LiuRo}), that $Y=A(\cdot, X) + c_2 X$ and $G=B(X)$.
It\^o's formula for the square of the $H$-norm applied to
\eqref{eq:1app} yields
\begin{align*}
  &\frac12 \norm[\big]{e^{-c_2 \cdot}X_\lambda}^2 + 
    \int_0^\cdot e^{-2c_2s} \ip[\big]{\tilde{A}_\lambda X_\lambda(s))}%
    {X_\lambda(s)}\,ds
    - \frac12 \int_0^\cdot e^{-c_2s}
    \norm[\big]{B_\lambda(X_\lambda(s))}^2_{\cL^2(U;H)}\,ds\\
  &\hspace{3em}= \frac12 \norm{X_0}^2
  + \int_0^\cdot e^{-c_2s}X_\lambda(s) B_\lambda(X_\lambda(s))\,dW(s).
\end{align*}
Since the stochastic integral on the right-hand side is a martingale,
as seen above, one has, applying inequality \eqref{eq:yo} and taking
expectations on both sides,
\begin{align*}
  &\frac12 \E\norm[\big]{e^{-c_2 \cdot}X_\lambda}^2
    - \frac12 \E\norm{X_0}^2\\
  &\hspace{3em} \leq \E\int_0^\cdot e^{-2c_2s} \Bigl(%
    -\ip[\big]{\tilde{A} J_\lambda X_\lambda(s)}{J_\lambda X_\lambda(s)}
    + \frac12 \norm[\big]{B(J_\lambda X_\lambda(s))}^2_{\cL^2(U;H)}%
    \Bigr)\,ds.
\end{align*}
Let $\varphi \in L^2(\Omega;C([0,T];H)) \cap
L^p(\Omega;L^p(0,T;V))$. Thanks to assumptions {\bf (III)--(IV)}, one
has that $\tilde A\varphi \in L^q(\Omega;L^q(0,T;V'))$ and
$B(\varphi) \in L^2(\Omega;L^2(0,T;\cL^2(U;H)))$, so that the term
within parentheses on right-hand side can be rewritten as
\begin{align*}
  &- \ip[\big]{\tilde{A}J_\lambda X_\lambda - \tilde{A}\varphi}
    {J_\lambda X_\lambda - \varphi}
    + \frac12 \norm[\big]{B(J_\lambda X_\lambda) - B(\varphi)}^2_{\cL^2(U;H)}\\
  &\hspace{3em} - \ip[\big]{\tilde{A}\varphi}{J_\lambda X_\lambda}
  - \ip[\big]{\tilde{A}J_\lambda X_\lambda - \tilde{A}\varphi}{\varphi}\\
  &\hspace{3em} - \frac12 \norm[\big]{B(\varphi)}^2_{\cL^2(U;H)}
  + \ip[\big]{B(J_\lambda X_\lambda)}{B(\varphi)}_{\cL^2(U;H)},
\end{align*}
where the sum of the two terms in the first row is negative by
assumption {\bf (III)}. One is thus left with
\begin{align*}
  &\frac12 \E\norm[\big]{e^{-c_2 \cdot} X_\lambda}^2
    - \frac12 \E\norm{X_0}^2\\
  &\hspace{3em}\leq \E\int_0^\cdot e^{-2c_2s} \Bigl(
  - \ip[\big]{\tilde{A}\varphi(s)}{J_\lambda X_\lambda(s)}
  - \ip[\big]{\tilde{A}J_\lambda X_\lambda(s) - \tilde{A}\varphi(s)}
  {\varphi(s)}\\
  &\hspace{5em} - \frac12 \norm[\big]{B(\varphi(s))}^2_{\cL^2(U;H)}
    + \ip[\big]{B(J_\lambda X_\lambda(s))}{B(\varphi(s))}_{\cL^2(U;H)}
    \Bigr)\,ds,
\end{align*}
from which it follows that, for every nonnegative
$\psi \in L^\infty(0,T)$,
\begin{align*}
  &\frac12 \E\int_0^T \psi(t) \bigl(e^{-2c_2t}
  \norm{X_\lambda(t)}^2 - \norm{X_0}^2 \bigr)\,dt\\
  &\hspace{3em} \leq
  \E\int_0^T \psi(t) \biggl(\int_0^t e^{-2c_2s} \Bigl(%
  - \ip[\big]{\tilde{A}\varphi(s)}{J_\lambda X_\lambda(s)}
    - \ip[\big]{\tilde{A}J_\lambda X_\lambda(s)%
    - \tilde{A}\varphi(s)}{\varphi(s)}\\
  &\hspace{5em} - \frac12 \norm[\big]{B(\varphi(s))}^2_{\cL^2(U;H)}
    + \ip[\big]{B(J_\lambda X_\lambda(s))}{B(\varphi(s))}_{\cL^2(U;H)}
    \Bigr)\,ds\biggr)\,dt.
\end{align*}
By the weak lower semicontinuity of the norm in $L^2(\Omega;L^2(0,T;H))$,
one has
\[
  \E\int_0^T \psi(t) \bigl(
  e^{-2c_2t} \norm{X(t)}^2 - \norm{X_0}^2 \bigr)\,dt
  \leq \liminf_{\lambda \to 0} \E\int_0^T \psi(t) \bigl(
  e^{-2c_2t} \norm{X_\lambda(t)}^2 - \norm{X_0}^2 \bigr)\,dt.
\]
Moreover, since $B(\varphi) \in L^2(\Omega;L^2(0,T;\cL^2(U;H)))$,
the weak convergence results proved above yield
\begin{align*}
  &\lim_{\lambda \to 0}
  \E\int_0^T \psi(t) \biggl( \int_0^t e^{-2c_2s} \Bigl(
  - \ip[\big]{\tilde{A}\varphi(s)}{J_\lambda X_\lambda(s)}
  - \ip[\big]{\tilde AJ_\lambda X_\lambda(s) - \tilde{A}\varphi(s)}{\varphi(s)}\\
  &\hspace{5em} - \frac12 \norm[\big]{B(\varphi(s))}^2_{\cL^2(U;H)}
    + \ip[\big]{B(J_\lambda X_\lambda(s))}{B(\varphi(s))}_{\cL^2(U;H)} \Bigr)\,ds
    \biggr)\,dt\\
  &= \E\int_0^T \psi(t) \biggl( \int_0^t e^{-2c_2s} \Bigl(
    - \ip[\big]{\tilde{A}(\varphi(s))}{X(s)}
    - \ip[\big]{Y(s)-\tilde{A}(\varphi(s))}{\varphi(s)}\\
  &\hspace{5em} - \frac12 \norm[\big]{B(\varphi(s))}^2_{\cL^2(U;H)}
  + \ip[\big]{G(s)}{B(\varphi(s))}_{\cL^2(U;H)} \Bigr)\,ds \biggr)\,dt.
\end{align*}
We then deduce that
\begin{align*}
  &\frac12 \E\int_0^T \psi(t) \bigl(e^{-2c_2t} \norm{X(t)}^2 - \norm{X_0}^2
    \bigr)\,dt\\
  &\hspace{3em} \leq \E\int_0^T \psi(t) \biggl( \int_0^te^{-2c_2s}
    \Bigl( - \ip[\big]{\tilde{A}\varphi(s)}{X(s)}
    - \ip[\big]{Y(s)-\tilde{A}\varphi(s)}{\varphi(s)}\\
  &\hspace{5em} - \frac12 \norm[\big]{B(\varphi(s))}^2_{\cL^2(U;H)}
  + \ip[\big]{G(s)}{B(\varphi(s))}_{\cL^2(U;H)} \Bigr)\,ds \biggr)\,dt.
\end{align*}
It\^o's formula (see \cite{KR-spde} or \cite[{\S}4]{LiuRo}) applied to
the limit equation \eqref{eq:lim} implies that there exists a
modification of $X$, denoted by the same symbol for simplicity, such
that $X \in L^2(\Omega;C([0,T];H))$ and
\begin{align*}
  &\frac12 e^{-2c_2 \cdot} \norm{X}^2
    + \int_0^\cdot e^{-2c_2s} \ip[\big]{Y(s)}{X(s)}\,ds
  - \frac12 \int_0^\cdot e^{-2c_2s}\norm[\big]{G(s)}^2_{\cL^2(U;H)}\,ds\\
  &\hspace{3em} = \frac12 \norm{X_0}^2
    + \int_0^\cdot e^{-2c_2s}X(s) G(s)\,dW(s).
\end{align*}
Substituting this identity on the left-hand side of the last
inequality and rearranging the terms yields
\[
  \E\int_0^T \psi(t) \int_0^t e^{-2c_2s} \Bigl(%
  \ip[\big]{Y(s)-\tilde{A}\varphi(s)}{\varphi(s)-X(s)}%
  + \frac12 \norm[\big]{G(s)-B(\varphi(s))}_{\cL^2(U;H)}^2%
  \Bigr)\,ds\,dt \leq 0.
\]
Choosing $\varphi=X$ immediately yields $G=B(X)$, hence also, in
particular,
\[
  \E\int_0^T \psi(t) \int_0^te^{-2c_2s}
  \ip[\big]{Y(s)-\tilde{A}\varphi(s)}{\varphi(s)-X(s)}\,ds\,dt
  \leq 0.
\]
Let $\delta \in \erre_+$, $v \in V$, and
$\bar{\varphi} \in L^\infty(\Omega \times [0,T])$. Choosing
$\varphi=X+\delta\bar{\varphi} v$ and taking the limit as
$\delta \to 0$ yields
\[
  \E\int_0^T\psi(t)\int_0^te^{-2c_2s} \bar{\varphi}(s)
  \ip[\big]{Y(s)-\tilde{A}\varphi(s)}{v}\,ds\,dt \leq 0.
\]
By a classical localisation argument, recalling that $\psi$ and
$\bar{\varphi}$ have been chosen arbitrarily, one has
\[
 \ip[\big]{Y-\tilde{A}\varphi}{v} \leq 0 \qquad \forall v \in V
\]
a.e. in $\Omega \times [0,T]$. Since $v$ has also been chosen
arbitrarily in $V$, it follows that
\[
  \ip[\big]{Y-\tilde{A}\varphi}{X - \varphi} \geq 0
\]
a.e. in $\Omega \times [0,T]$. Since $\tilde{A}$ is maximal monotone,
this implies that $Y=\tilde{A}X$ a.e. in $\Omega \times [0,T]$. We
have thus shown that $X$ is a strong solution to equation
\eqref{eq:1}.

In order to prove the Lipschitz continuous dependence of the solution
with respect to the initial datum (from which uniqueness follows), let
$X_1$, $X_2$ be strong solutions to \eqref{eq:1}, in the sense of
Theorem~\ref{th:wp}, with initial data $X_0^1$ and $X_0^2$,
respectively. It\^o formula (as in \cite{KR-spde} or \cite{LiuRo})
yields
\begin{align*}
  &\frac12 \norm[\big]{(X_1-X_2)}^2_H\\
  &\hspace{1em} + \int_0^\cdot \Bigl(\ip[\big]{AX_1-AX_2}{X_1-X_2}
  - \frac12 \norm[\big]{B(X_1)-B(X_2)}_{\cL^2(U;H)}^2\Bigr)(s)\,ds\\
  &\hspace{3em} = \frac12 \norm[\big]{X_0^1-X_0^2}^2
  + \int_0^\cdot \bigl((X_1-X_2)(B(X_1)-B(X_2))\bigr)(s)\,dW(s),
\end{align*}
where the stochastic integral on the right-hand side is a martingale
because $X_1$, $X_2 \in L^2(\Omega;C([0,T];H))$ and $B(X_1)$,
$B(X_2) \in L^2(\Omega;L^2(0,T;\cL^2(U;H)))$. Therefore, taking
expectations on both sides and using assumption \textbf{(III)},
\[
  \E\norm[\big]{(X_1-X_2)(t)}^2 \leq
  \E\norm[\big]{X_0^1-X_0^2}^2
  + 2c_2\int_0^t \E\norm[\big]{(X_1-X_2)(s)}^2\,ds \qquad\forall t \in [0,T],
\]
from which the conclusion follows thanks to Gronwall's inequality. The
proof of Theorem~\ref{th:wp} is thus complete.

\bigskip\par\noindent
\textbf{Acknowledgements.} CM gratefully acknowledges the hospitality
of the Interdisziplin\"ares Zentrum f\"ur Komplexe Systeme,
Universit\"at Bonn, Germany. LS and US are partially funded by the
FWF-DFG grant I\,4354, the FWF grants F\,65, P\,32788, and M 2876, and
the Vienna Science and Technology Fund (WWTF) through Project
MA14-009.


\bibliographystyle{amsplain}

\def\polhk#1{\setbox0=\hbox{#1}{\ooalign{\hidewidth
  \lower1.5ex\hbox{`}\hidewidth\crcr\unhbox0}}}
\providecommand{\bysame}{\leavevmode\hbox to3em{\hrulefill}\thinspace}
\providecommand{\MR}{\relax\ifhmode\unskip\space\fi MR }
\providecommand{\MRhref}[2]{%
  \href{http://www.ams.org/mathscinet-getitem?mr=#1}{#2}
}
\providecommand{\href}[2]{#2}

\end{document}